\documentclass[a4paper,12pt,leqno]{article}
\usepackage[utf8]{inputenc}
\usepackage{amssymb,amsthm}
\usepackage{scalerel}
\usepackage{hyperref,color}
\usepackage{url}
\hypersetup{
	colorlinks=true,
	linkcolor=blue,
	filecolor=magenta,
	urlcolor=cyan,}
\usepackage[left=2.3cm,top=3cm,right=2.3cm,bottom=3cm,bindingoffset=0.5cm]{geometry}
 \newcommand\blfootnote[1]{%
	\begingroup
	\renewcommand\thefootnote{}\footnote{#1}%
	\addtocounter{footnote}{-1}%
	\endgroup
}

\usepackage{titlesec}
\titlelabel{\thetitle.\quad}
\titleformat*{\section}{\large\bf}
\titleformat*{\subsection}{\normalsize\bf}
\title{DOUBLE TAILS OF MULTIPLE ZETA VALUES}
\author{P. Akhilesh}
\title{DOUBLE TAILS OF MULTIPLE ZETA VALUES}
\newcommand{\Addresses}{{
		\bigskip
		\footnotesize
		P.~Akhilesh, \textsc{Harish-Chandra Research Institute, Jhunsi, Allahabad -211 019, India. \\ Institute of Mathematical Sciences, Chennai ,  India}\par\nopagebreak 
		\textit{E-mail }: \texttt{akhi@imsc.res.in}
}}
\date{}
\font\pc=cmcsc10
\newtheorem{thm}{\pc Theorem }
\newtheorem{Def}{\pc Definition}
\newtheorem{prop}{\pc Proposition}
\newtheorem*{cor}{\pc Corollary}
\def\cvirg{\,\raise2pt\hbox{,}}

\begin{document}
\maketitle
   

\begin{abstract}
In this paper we introduce and study double tails of multiple zeta values. We show, in particular, 
that they satisfy certain recurrence relations and  deduce from this a generalization of Euler's classical 
formula  $\zeta(2)=3\sum_{m=1}^\infty m^{-2}\left({2m\atop m}\right)^{-1}$ to all multiple zeta values, as well as  a new and very efficient algorithm for computing these values.

\end{abstract}

\section{Introduction}
\indent 

Throughout the paper, $\bf N$ denotes the set of non-negative integers. A finite sequence ${\bf a}=(a_1,\ldots,a_r)$ of positive integers  is called {\it a composition}. The integer $r$ is called {\it the depth} of $\bf a$ and the integer $k=a_1+\ldots+a_r$ {\it the weight} of $\bf a$. The composition~$\bf a$ is said to be {\it admissible} if either
$r\geqslant 1$ and $a_1\geqslant 2$, or $\bf a$ is the empty composition denoted $\varnothing$.

To each admissible composition ${\bf a}=(a_1,\ldots,a_r)$, one associates a real number~$\zeta(\bf a)$. It is defined by the convergent series
\begin{equation}
\label{E1}
\zeta({\bf a})= \sum_{n_1>\ldots>n_r>0}n_1^{-a_1}\ldots n_r^{-a_r}.
\end{equation}
when $r\geqslant 1$, and by $\zeta(\varnothing)=1$ when $r=0$. These numbers are called {\it multiple zeta values} or
{\it Euler-Zagier numbers}.

A {\it binary word}  is by definition a word $w$ constructed on the alphabet $\{0,1\}$. Its letters are  called {\it bits}. The number of bits of $w$ is called {\it the weight} of $w$ and denoted by~$|w|$. The number of bits of $w$ equal to $1$ is called {\it the depth} of $w$.  To any composition ${\bf a}=(a_1,\ldots,a_r)$, one associates the binary word
\begin{equation}
\label{E2}
{\bf w({\bf{a}})}= \{0\}_{a_1-1}1\ldots \{0\}_{a_r-1}1
\end{equation}
where for each integer $u\geqslant 0$, $\{0\}_u$ denotes the binary word consisting of $u$ bits equal to~$0$, and where $\bf w(\bf a)$ is the empty word if $\bf a$ is the empty composition.
The weight of $\bf w(\bf a)$ is equal to the weight of $\bf a$ and its depth to the depth of $\bf a$. 

We shall denote by ${\rm W}$ the set of binary words. When $\varepsilon,\varepsilon'\in\{0,1\}$,
${}_\varepsilon\!{\rm W}$ and ${\rm W}\!_{\varepsilon'}$ denote the sets of binary words starting by $\varepsilon$ and ending by $\varepsilon'$ respectively,
and ${}_\varepsilon\!{\rm W}\!_{\varepsilon'}$ their intersection.

The map $\bf w$ is a bijection from the set of compositions onto the set of binary words not ending by $0$. Non empty compositions correspond to words in ${\rm W}\!_1$, and non empty admissible compositions to words 
in ${}_0\!{\rm W}\!_1$. Therefore a binary word will be called {\it admissible} if either it belongs to 
${}_0\!{\rm W}\!_1$, or it is empty.

Maxim Kontsevich has discovered that for each  admissible composition $\bf a$, the multiple zeta value $\zeta(\bf a)$ can be written as an iterated integral. More precisely, if $w=\varepsilon_1\ldots\varepsilon_k$ denotes  the associated binary word $\bf w(\bf a)$, we have
\begin{equation}
\label{E3}
\zeta({\bf{a}})={\rm{It}}\int_0^1(\omega_{\varepsilon_1},\ldots,\omega_{\varepsilon_k})= \int_{1>t_1>\ldots>t_k>0} f_{\varepsilon _1}(t_1)\ldots f_{\varepsilon _k}(t_k)dt_1\ldots dt_k
\end{equation}
where $\omega_i=f_i(t)dt$, with  $f_0(t)={1\over t}$ and $f_1(t)={1\over 1-t}$. We therefore often
write this number $\zeta(w)$ instead of $\zeta(\bf a)$.

Let $w=\varepsilon_1\ldots\varepsilon_k$ be a binary word. Its {\it dual word} is defined to be $\overline{w}=\overline{\varepsilon}_k\ldots\overline{\varepsilon}_1$, where $\overline{0}=1$ and $\overline{1}=0$. When $w$ is admissible, so is $\overline{w}$. We can therefore define
the {\it dual composition} of an admissible composition $\bf a$ to be the admissible composition $\overline{\bf a}$ such that $\bf w(\overline{\bf a})$ is dual to $\bf w(\bf a)$. When $\bf a$ has weight $k$ and depth $r$, $\overline{\bf a}$ has weight $k$ and depth~$k-r$. 

By the change of variables $t_i\mapsto 1-t_{k+1-i}$ in the integral (\ref{E3}), one gets the following {\it duality relation} : for any admissible 
composition $\bf a$, we have
\begin{equation}
\label{E4}
\zeta(\bf a)=\zeta(\overline{\bf a})
\end{equation}
When $\bf a$ is a non empty admissible composition, we can define for each integer $n\geqslant 0$ {\it the $n$-tail} of the series (\ref{E1})
to be the sum of the series 
\begin{equation}
\label{E5}
\sum_{n_1>\ldots>n_r>n}n_1^{-a_1}\ldots n_r^{-a_r}.
\end{equation}

We remark in Section \ref{S1} that this $n$-tail can be written as the iterated integral
\begin{equation}
\label{E6}
{\rm It}\int_0^1(\omega_{\varepsilon_1},\ldots,t^n\omega_{\varepsilon_k})=\int_{1>t_1>\ldots>t_k>0}f_{\varepsilon _1}(t_1)\ldots
f_{\varepsilon _k}(t_k)t_k^ndt_1\ldots dt_k
\end{equation}
where $\varepsilon_1\ldots\varepsilon_k$ is the binary word $\bf w(\bf a)$. 
This observation  motivates the following definition:\medskip
\begin{Def}---  \label{D1} When $\bf a$ is a non empty admissible composition, we define for  $m$ and $n$ in $\bf N$  the 
{ $(m,n)$-double tail} $\zeta({\bf{a}})_{m,n}$ of $\zeta(\bf a)$ as the iterated integral
\begin{eqnarray}
\label{E7}
\zeta({\bf{a}} )_{m,n}  &=&{\rm It} \int_0^1((1-t)^m\omega_{\varepsilon_1},\ldots,t^n\omega_{\varepsilon_k})  \nonumber \\&=&\int_{1>t_1>\ldots>t_k>0}(1-t_1)^m f_{\varepsilon _1}(t_1)\ldots 
f_{\varepsilon_k}(t_k) t_k^n dt_1\ldots  dt_k \,,
\end{eqnarray}
where $\varepsilon_1\ldots\varepsilon_k$ is the binary word $\bf w(\bf a)$. \medskip
\end{Def}
Our purpose in this paper is to study  these double tails. We first establish a series expression for them in Section \ref{S1} :\medskip

\begin{thm}--- \label{T1} Let ${\bf a}=(a_1,\ldots,a_r)$ be a non empty admissible composition. For  all $m$ and  $n$ in~$\bf N$, the 
$(m,n)$-double tail of $\zeta(\bf a)$ is given by the convergent series
\begin{equation}
\label{E8}
\zeta({\bf a})_{m,n} =\sum_{n_1>\ldots>n_r>n}\left({n_1+m\atop m}\right)^{-1}n_1^{-a_1}\ldots n_r^{-a_r}.
\end{equation}
\end{thm}
By applying the change of variables $t_i\mapsto 1-t_{k+1-i}$ in the  integral (\ref{E7}), we get in Section~\ref{S2} the following duality relation for double tails :\medskip
\begin{thm}---
\label{T2}  Let ${\bf a}$ be a non empty admissible composition and $\overline{{\bf a}}$ denote its dual composition. For any
 $m$ and $n$ in $\bf N$, we have 
\begin{equation}
 \label{E9}
\zeta({\bf a})_{m,n}=\zeta(\overline{{\bf a}})_{n,m}
\end{equation}
 \end{thm}

Note that $\zeta({\bf a})_{0,n}$ is nothing but the usual $n$-tail of $\zeta(\bf a)$. Formula (\ref{E9}) tells us that it is equal $\zeta(\overline{\bf a})_{n,0}$.
This equality is in fact the main theorem of a recent paper   by J.~M.~Borwein and O-Yeat Chan
(\cite{2}, th. 14),  for which Theorem \ref{T2} therefore provides a conceptually very simple  proof. \medskip

When ${\bf a}=(a_1,\ldots,a_r)$ is a non empty admissible composition, the $n$-tail $\zeta({\bf a})_{0,n}$ of $\zeta({\bf a})$ 
is only slowly decreasing to $0$ when $n$ grows. More precisely, as we shall see in section \ref{S3}, $\zeta({\bf a})_{0,n}$ is equivalent to 
\begin{equation}
\label{E10}
 {n^{r-(a_1+\ldots+a_r)}\over (a_1-1)(a_1+a_2-2)\ldots(a_1+\ldots+a_r-r)}
\end{equation}
when $n$ tends to $+\infty$. The double tails of $\zeta(\bf a)$ are generally much smaller, as is shown by the following theorem, also proved in Section \ref{S3} :\medskip

\begin{thm}\label{T3}---  Let $\bf a$ be a non empty admissible composition. For all $m$ and $n$ in $\bf N$, we have

\begin{equation}
\label{E11}
\zeta({\bf a})_{m,n}\leqslant {m^mn^n\over (m+n)^{m+n}}\zeta(\bf a),
\end{equation}
and $\zeta({\bf a})\leqslant{\pi^2\over 6}$. \it{We have in particular}
\begin{equation}
\label{E12}
\zeta({\bf a})_{n,n}\leqslant 2^{-2n}\zeta({\bf a})\leqslant 2^{-2n}{\pi^2\over 6}.
\end{equation}

\end{thm}

In the following sections of the paper, we study the recurrence relations satisfied by  double tails of multiple zeta values. In order to state them
in a nice and simple form, it is convenient to use binary words rather than compositions and to slightly extend Definition~\ref{D1} as follows.\medskip

\begin{Def}\label{D2}---  Let $w=\varepsilon_1\ldots\varepsilon_k$ be a binary word and let $m,n\in\bf N$. Assume $m\geqslant 1$ when $w\in{}_1\!{\rm W}$,  and $n\geqslant 1$ when $w\in{\rm W}\!_{0}$. We define a real number $\zeta(w)_{m,n}$ by the convergent iterated integral
\begin{equation}
\label{E13}
\zeta(w)_{m,n}={\rm It}\int_0^1((1-t)^m\omega_{\varepsilon_1},\ldots,t^n\omega_{\varepsilon_k}),
\end{equation}
when $k\geqslant 2$, and in the remaining cases by
\begin{equation}
\label{E14}
\zeta(0)_{m,n}=\int_0^1(1-t)^{m}t^{n}{dt\over t}={m!(n-1)!\over (m+n)!}\,,\quad\;
\end{equation}
\begin{equation}
\label{E15}
\zeta(1)_{m,n}=\int_0^1(1-t)^{m}t^{n}{dt\over 1-t}={(m-1)!n!\over (m+n)!}\,,
\end{equation}
\begin{equation}
\label{E16}
\zeta(\varnothing)_{m,n}={m!\,n!\over (m+n)!}\cdot \qquad\qquad \qquad \qquad\; \qquad
\end{equation}
\end{Def}
The hypotheses $m\geqslant 1$ when $w\in{}_1\!{\rm W}$  and $n\geqslant 1$ when $w\in{\rm W}\!_{0}$ are needed in order to ensure that the integrals involved converge.
Of course, when $w$ is the binary word ${\bf w({\bf a})}$ associated to a non empty admissible composition $\bf a$, the number
$\zeta(w)_{m,n}$ is  the same as the double tail $\zeta({\bf a})_{m,n}$. Moreover, whenever $\zeta(w)_{m,n}$ is defined, so is 
$\zeta(\overline{w})_{n,m}$ and we have the duality relation
\begin{equation}
\label{E17}
\zeta(w)_{m,n}=\zeta(\overline{w})_{n,m}.
\end{equation}

We now state the  recurrence relations satisfied by the numbers $\zeta(w)_{m,n}$ :\medskip

\begin{thm}
\label{T4}---  Let $w$ be a binary word and and let $m,n\in\bf N$.
\\$a)$ Assume $n\geqslant 1$. Then we have
\begin{equation}
\label{E18}
\cases{\zeta(w0)_{m,n}=n^{-1}\zeta(w)_{m,n}&if $m\geqslant 1$ or $w\notin{}_1\!{\rm W}$,\cr
\zeta(w1)_{m,n-1}=\zeta(w1)_{m,n}+n^{-1}\zeta(w)_{m,n}&if $m\geqslant 1$ or $w\in{}_0\!{\rm W}$.\cr}
\end{equation}
$b)$ Assume $m\geqslant 1$.  Then we have
\begin{equation}
\label{E19}
\cases{\zeta(1w)_{m,n}=m^{-1}\zeta(w)_{m,n}&if $n\geqslant 1$ or $w\notin{\rm W}\!_0$,\cr
\zeta(0w)_{m-1,n}=\zeta(0w)_{m,n}+m^{-1}\zeta(w)_{m,n}&if $n\geqslant 1$ or $w\in{\rm W}\!_1$.\cr}
\end{equation}
\par \medskip
\end{thm}
 \indent{\it{Remark}}---  These recurrence relations are slightly cumbersome to state, because of the hypotheses  which we have to impose on $m$ and $n$ in order  to ensure that all terms involved  are well defined. One could however also extend Definition \ref{D2} further and define~$\zeta(w)_{m,n}$ for all binary words $w$ and all integers $m$ and $n$ in $\bf N$, by allowing in formulas (\ref{E13}), (\ref{E14}), (\ref{E15}) regularized iterated integrals instead of only convergent ones. In that case for example, $\zeta(1)_{0,n}$ would be equal to $-\!\sum_{1\leqslant k\leqslant n}{1\over k}$ for all $n\in\bf N$ and 
formula (\ref{E18}) (resp. formula  (\ref{E19})) would hold with no restriction on $m$ (resp. on $n$).\medskip

Theorem \ref{T4} will be proved in Section \ref{S4}, together with some of its consequences. Worth mentioning is the following
one :\medskip

\begin{thm}\label{T5}--- {\it Let $w$ be an non empty admissible binary word. There exists a unique triple $(v,a,b)$, where $v$ is an admissible binary word, empty or not, and $a$, $b$ are positive integers, such that $w=0\{1\}_{b-1}v\{0\}_{a-1}1$. For each integer $n\geqslant 1$, we have
\begin{equation}
\label{E20}
\hspace{0.5cm}
\zeta(w)_{n-1,n-1}=\zeta(w)_{n,n}+n^{-a}\zeta(w^{\rm init})_{n,n}+n^{-b}\zeta(w^{\rm fin})_{n,n}+n^{-a-b}\zeta(w^{\rm mid})_{n,n}
\end{equation}
where $w^{\rm init}=0\{1\}_{b-1}v$, $w^{\rm fin}=v\{0\}_{a-1}1$ and $w^{\rm mid}=v$.}\medskip

\end{thm}

The words $w^{\rm init}$, $w^{\rm fin}$ and $w^{\rm mid}$ will be called the {\it initial part}, the {\it final part} and {\it the middle part} of the word $w$. Recurrence relations (\ref{E20}), together with upper bounds (\ref{E12}) lead to a very simple and  fast algorithm to compute multiple zeta values,
as will be explained in Section \ref{S5}. This algorithm is especially efficient when one is interested in getting simultaneously all multiple zeta values
up to a given weight $k$. \medskip

When applied to the binary word $01$, corresponding to the composition $(2)$, Theorem~\ref{T5} yields
$\zeta(2)_{n-1,n-1}=\zeta(2)_{n,n}+3n^{-2}\big({2n\atop n}\big)^{-1}$ for $n\geqslant 1$, hence by  induction \\${\zeta(2)=\zeta(2)_{n,n}+3\sum_{1\leqslant m\leqslant n} m^{-2}\big({2m\atop m}\big)^{-1}}$ for $n\geqslant 0$. Letting $n$ go to~$+\infty$, we get the following formula, discovered by Euler
\begin{equation}
\label{E21}
\zeta(2)=3\sum_{m=1}^\infty m^{-2}\left({2m\atop m}\right)^{-1}.
\end{equation}
Other formulas of the same type are also easily deduced from Theorem \ref{T5}, such as for example
\begin{equation}
\label{E22}
\zeta(4)={36\over17}\sum_{m=1}^{\infty} m^{-4}\left({2m\atop m}\right)^{-1},
\end{equation}
which seems to appear  in print  for the first time as an exercise in a book by L. Comtet in 1974 (\cite{3}, p. 89), and 
was then experimentally rediscovered and finally  proved by A.~van~der~Poorten in 1979 (\cite{5}, p. 203); see also  P. Bala (\cite{1}, p. 15). These exemples
are treated in Section \ref{S6}.

In Section \ref{S7}, we extend these types of formulas to all multiple zeta values. To state our result, we shall need some further notations. For any non empty composition 
${{\bf a}=(a_1,\ldots,a_r)}$ (admissible or not) and any integer $m\geqslant 1$, we define a real number $\varphi_m(\bf a)$ by the finite sum
\begin{equation}
\label{E23}
\varphi_m({\bf a})=m^{-a_1}\sum_{m>n_2>\ldots>n_r>0}n_2^{-a_2}\ldots n_r^{-a_r},
\end{equation}
where the right hand side is considered to be equal to $m^{-a_1}$ when $r=1$. For any pair $(\varepsilon,\varepsilon ')$ of bits,
we define an integer $\lambda(\varepsilon,\varepsilon')$ as follows
\begin{equation}
\label{E24}
\lambda(\varepsilon,\varepsilon')=\cases{1&if $(\varepsilon,\varepsilon ')$ is equal to $(1,0)$,\cr
2&if $(\varepsilon,\varepsilon ')$ is equal to $(0,0)$ or $(1,1)$,\cr
3&if $(\varepsilon,\varepsilon ')$ is equal to $(0,1)$,\cr}
\end{equation}
With these notations :\medskip

\begin{thm}--- \label{T6} Let $\bf a$ be a non empty admissible composition. Let $\varepsilon_1\ldots\varepsilon_k$ denote the corresponding binary word
$\bf w(\bf a)$. For any index i such that $1\leqslant i\leqslant  k-1$, let ${\bf a}_i$ and ${\bf b}_i$ denote the compositions corresponding to the binary words
$\varepsilon_{i+1}\ldots\varepsilon_k$ and $\overline{\varepsilon_i}\ldots\overline{\varepsilon_1}$ respectively. We then have
\begin{equation}
\label{E25}
\zeta({\bf a})=\sum_{m=1}^{\infty}\psi_m({\bf a})\left({2m\atop m}\right)^{-1}
\end{equation}
where for each $m\geqslant 1$
\begin{equation} 
\label{E26}
\psi_m({\bf a})=\sum_{i=1}^{k-1}\lambda(\varepsilon_i,\varepsilon_{i+1})\varphi_m({\bf a}_i)\varphi_m({\bf b}_i).
\end{equation}
\end{thm}
Theorem \ref{T6} yields a second algorithm to compute multiple zeta values. It is faster than the first one when one is only interested in computing a single
multiple zeta value. It is also faster than the standard algorithm, based on evaluations
of multiple polylogarithms at ${1\over 2}$, which has been for example implemented by J.~Borwein, P.~Lisonek, P.~Irvine and C.~Chan on their website EZ-Face \footnote{\hskip 0.1cm {\url{ http://wayback.cecm.sfu.ca/projects/EZFace/Java/index.html}\par}}.
\noindent More details on this topic will
be given in Section~\ref{S8}.

As Henri Cohen pointed us in a private communication~\cite{cohen}, our algorithms can be extended to  compute values of multiple polylogarithms (in many variables). We are most grateful to him for providing us with this possible extension of our work and we hope to address this question in depth in the future.
\bigskip
\section{Series expression for double tails of multiple zeta values}\label{S1}\smallskip

We shall first prove that the iterated integrals considered in Definition \ref{D2} converge:\medskip

\begin{prop}--- \label{P1} Let $w=\varepsilon_1\ldots\varepsilon_k$ be a binary word of weight $k\geqslant2$, and let $m,n\in\bf N$. Assume $m\geqslant1$ when $\varepsilon_1=1$ and $n\geqslant 1$ when $\varepsilon_k=0$. The iterated integral (\rm{\ref{E13}}) defining
$\zeta(w)_{m,n}$ is convergent.
\end{prop}

 By definition, this iterated integral can be expressed as
\begin{equation}
\label{E27}
\int_{1>t_1>\ldots>t_k>0}(1-t_1)^mf_{\varepsilon_1}(t_1)\ldots f_{\varepsilon_k}(t_k)t_k^n dt_1\ldots dt_k,
\end{equation}
were $f_0(t)={1\over t}$ and $f_1(t)={1\over 1-t}$. Our hypotheses on $m$ and $n$ imply that, on the domain of integration, we have
\begin{equation}
\label{E28}
(1-t_1)^mf_{\varepsilon_1}(t_1)\leqslant f_0(t_1),\qquad\qquad f_{\varepsilon_k}(t_k)t_k^n\leqslant f_1(t_k)\,,
\end{equation}
hence the integral (\ref{E27}) is bounded above by the integral
\begin{equation}
\label{E29}
\int_{1>t_1>\ldots>t_k>0}f_{0}(t_1)f_{\varepsilon_2}(t_2)\ldots f_{\varepsilon_{k-1}}(t_{k-1})f_{1}(t_k) dt_1\ldots dt_k,
\end{equation}
which by (\ref{E3}) converges to $\zeta(w')$, where $w'$ is the admissible binary word $0\varepsilon_2\ldots\varepsilon_{k-1}1$. \medskip
 
 We now state some elementary properties of the numbers $\zeta(w)_{m,n}$. Note that  their expression (\ref{E27}) 
 is valid even when the weight of $w$ is $1$, as formulas (\ref{E14}) and (\ref{E15}) show.
  \medskip

\begin{prop}---  \label{P2}Let $w$ be a binary word and let $m\in\bf N$.
\par
$a)$  Assume $m\geqslant1$ if $w\in{}_1\!{\rm W}$. Then $\zeta(w0)_{m,n}=n^{-1}\zeta(w)_{m,n}$ for  $n\geqslant 1$.\par
$b)$ Assume $m\geqslant1$ or $w\in{}_0\!{\rm W}$. Then $\zeta(w1)_{m,n}=\sum_{\ell> n} \zeta(w0)_{m,\ell}$ for  $n\geqslant 0$.
\end{prop}

Formulas (\ref{E14}) and (\ref{E16}) imply assertion $a)$ when $w$ is the empty binary word. When $w=\varepsilon_1\ldots\varepsilon_k$ is non empty, we  express  $\zeta(w)_{m,n}$ and $\zeta(w0)_{m,n}$ as integrals by (\ref{E27}).  Assertion $a)$ then follows form the identity
$\int_0^{t_k}t_{k+1}^n{dt_{k+1}\over t_{k+1}}={t_{k}^n\over n}$, valid for all $t_k\in]0,1[$ and~$n\geqslant 1$. Similarly, assertion $b)$
follows from the identity \smash{${t_{k+1}^n\over 1-t_{k+1}}=\sum_{\ell> n} t_{k+1}^{\ell-1}$},  valid for all $t_{k+1}\in]0,1[$ and $n\geqslant 0$.\medskip

\begin{cor}---  Let $w$ be a binary word and let $m\geqslant 0$, $n\geqslant 0$,
$a\geqslant 1$ be integers. Assume that one of the three following hypotheses holds : $m\geqslant 1$; $w\in{}_0\!{\rm W}$; $w=\varnothing$ and $a\geqslant 2$. Then
\begin{equation}
\label{E30}
\zeta(w\{0\}_{a-1}1)_{m,n}=\sum_{\ell> n}\ell^{-a}\zeta(w)_{m,\ell}\,.
\end{equation}
\end{cor}

Indeed, we have $\zeta(w\{0\}_{a-1}1)_{m,n}=\sum_{\ell> n}\zeta(w\{0\}_{a})_{m,\ell}$ by
Proposition \ref{P2}, $b)$, \\and~$\zeta(w\{0\}_{a})_{m,\ell}=\ell^{-a}\zeta(w)_{m,\ell}$ for $\ell>0$ by 
Proposition \ref{P2}, $a)$.\medskip

We now prove Theorem \ref{T1}. We have to show that, for any non empty admissible composition ${\bf a}=(a_1,\ldots,a_r)$ and  all $m$ and  $n$ in $\bf N$, the 
$(m,n)$-double tail $\zeta({\bf a})_{m,n}$ of $\zeta(\bf a)$ is given by the convergent series
\begin{equation}
\label{E31}
\sum_{n_1>\ldots>n_r>n}\left({n_1+m\atop m}\right)^{-1}n_1^{-a_1}\ldots n_r^{-a_r}.
\end{equation}
The integral (\ref{E7}) defining   $\zeta({\bf a})_{m,n}$ converges by Proposition \ref{P1}. We proceed by induction on~$r$. When $r=1$, we apply Corollary  of Proposition \ref{P2} to
the empty binary word and get
$$\zeta(a_1)_{m,n}=\zeta(\{0\}_{a_1-1}1)_{m,n}=\sum_{n_1> n}n_1^{-a_1}\zeta(\varnothing)_{m,n_1}=
\sum_{n_1> n}\left({n_1+m\atop m}\right)^{-1}n_1^{-a_1}.$$
When $r\geqslant 2$,  we apply Corollary  of Proposition \ref{P2} to the binary word $\bf w({\bf a'})$ where ${\bf a'}=(a_1,\ldots,a_{r-1})$ and  
use the induction hypothesis. We get
\begin{eqnarray}
\zeta({\bf a})_{m,n}&=&\zeta({\bf w}({\bf a}))_{m,n}=\zeta({\bf w}({\bf a}')\{0\}_{a_r-1}1)_{m,n} 
 =\sum_{n_r> n}n_r^{-a_r}\zeta({\bf w}({\bf a}'))_{m,n_r}\quad
 \nonumber\\ 
&=&\sum_{n_r> n}n_r^{-a_r}\zeta({\bf a}')_{m,n_r} \quad 
 =
\sum_{n_1>\ldots>n_r>n}\left({n_1+m\atop m}\right)^{-1}n_1^{-a_1}\ldots n_r^{-a_r} \;\nonumber.
\end{eqnarray}

\section{Duality relation for double tails}\label{S2}\smallskip

Let $w=\varepsilon_1\ldots\varepsilon_k$ be a binary word and let $m,n\in\bf N$. Assume  $m\geqslant 1$ when $w\in{}_1\!{\rm W}$  and $n\geqslant 1$ when $w\in{\rm W}\!_{0}$. Then $\zeta(w)_{m,n}$ as well as $\zeta(\overline{w})_{n,m}$ are well defined. They satisfy the duality relation stated as formula (\ref{E17}) in the introduction :
$$\zeta(w)_{m,n}=\zeta(\overline{w})_{n,m}.$$
This equality follows indeed from formula (\ref{E16})  when $k=0$ and, when $k\geqslant 1$, from the change of variables 
$t'_i=1-t_{k+1-i}$ in the integral expression (\ref{E27}) of $\zeta(w)_{m,n}$.\medskip

When applying this  relation to the word  associated to a non empty admissible composition, we get  Theorem \ref{T2}.
\bigskip

\section{Upper bounds for double tails} \label{S3}

Let ${\bf a}=(a_1,\ldots,a_r)$ be a non empty admissible composition and $n$ be a non-negative integer. We recall that the $n$-tail of
$\zeta(\bf a)$ is by definition
 the sum of the series
$$\sum_{n_1>\ldots>n_r>n}n_1^{-a_1}\ldots n_r^{-a_r}.$$
It  is therefore equal to $\zeta({\bf a})_{0,n}$ by Theorem \ref{T1}, or equivalently to the iterated integral  (\ref{E6}).

We shall now prove that  this $n$-tail is equivalent to expression (\ref{E10}) when $n$ tends to~$+\infty$. 
We proceed by induction on $r$. When $r=1$, the $n$-tail of $\zeta(\bf a)$ is  $\sum_{n_1>n}n_1^{-a_1}$, and
is equivalent to the integral \smash{$\int_{n}^{+\infty}x^{-a_1}dx={n^{1-a_1}\over a_1-1}$}. Our assertion follows in this case.

Assume now  $r\geqslant 2$. Then $\zeta({\bf a})_{0,n}$ can be written as $\sum_{n_r>n}n_r^{-a_r}\zeta({\bf a'})_{0,n_r}$, where 
${{\bf a}'=(a_1,\ldots,a_{r-1})}$. Applying the induction hypothesis to $\bf a'$, we see that $\zeta({\bf a})_{0,n}$  is equivalent to
$$\sum_{n_r>n}{n_r^{r-1-a_1-\ldots-a_r}\over (a_1-1)(a_1+a_2-2)\ldots (a_1+\ldots+a_{r-1}-r+1)}$$
when $n$ tend to $+\infty$, and therefore  to (\ref{E10}) by the previous paragraph.\medskip

We shall now prove Theorem \ref{T3}, which gives upper bounds for the double tails of multiple zeta values. Let again 
${\bf a}=(a_1,\ldots,a_r)$ be a non empty admissible composition and let $m,n\in\bf N$. From the expressions
(\ref{E3}) and (\ref{E7}) of $\zeta(\bf a)$ and $\zeta({\bf a})_{m,n}$ as iterated integrals, one deduces that 
$$\zeta({\bf a})_{m,n}\leqslant c_{m,n}\zeta({\bf a}),$$
where $c_{m,n}$ is the upper bound of $(1-t_1)^mt_k^n$ when $1\geqslant t_1\geqslant t_k\geqslant 0$, or equivalently 
the maximum value of the function $f:t\mapsto (1-t)^mt^n$ on $[0,1]$. Since the logarithmic derivative of $f$ is ${m\over t-1}+{n\over t}$,
$f$ attains its maximum when $t={n\over m+n}$,  and therefore  $c_{m,n}$ is equal to ${m^m n^n\over (m+n)^{m+n}}$. This proves
formula (\ref{E11}) in Theorem \ref{T3}. Furthermore, it is known that $\zeta(a)\leqslant {\pi^2\over 6}$ (see for example \cite{4}, Theorem 1). Formula (\ref{E12}) follows by taking $m=n$.
\bigskip
\section{Recurrence relations for double tails}\label{S4}

The first formula of Theorem \ref{T4}, $a)$, has been proved in Proposition \ref{P2}, $a)$. The second one then follows from Proposition \ref{P2}, $b)$. Assertion $b)$ of Theorem \ref{T4} is deduced from assertion $a)$ by the duality relations 
(\ref{E17}).
 This completes the proof of the recurrence relations of Theorem~\ref{T4}. \medskip
 
We now derive some of their consequences :\medskip

\begin{prop}--  \label{P3}Let $v$ be a binary word and let $m\geqslant 0$, $n\geqslant 1$,
$a\geqslant 1$ be integers. Assume that one of the three following hypotheses holds : $m\geqslant 1$; $v\in{}_0\!{\rm W}$; $v=\varnothing$ and $a\geqslant 2$. Let $w$ denote the word $v\{0\}_{a-1}1$. Then
\begin{equation}
\zeta(w)_{m,n-1}=\zeta(w)_{m,n}+n^{-a}\zeta(v)_{m,n}.
\end{equation}
\end{prop}
\medskip
This  follows immediately from relations (\ref{E18}), or can as well be deduced from the Corollary of Proposition \ref{P2}. Similarly, one deduces from relations~(\ref{E19}), or by duality from Proposition~\ref{P3}:\medskip

\begin{prop}\label{P4}.-- \it Let $v$ be a binary word and let $m\geqslant 1$, $n\geqslant 0$,
$a\geqslant 1$ be integers. Assume that one of the three following hypotheses holds : $n\geqslant 1$; $v\in{\rm W}\!_1$; $v=\varnothing$ and $a\geqslant 2$. Let $w$ denote the word $0\{1\}_{b-1}v$. Then
\begin{equation}
\label{E33}
\zeta(w)_{m-1,n}=\zeta(w)_{m,n}+m^{-b}\zeta(v)_{m,n}.
\end{equation}

\end{prop}

We  now prove the following proposition, from which Theorem \ref{T5} is a particular case :\medskip

\begin{prop}\label{P5}.-- {\it Let $w$ be a non empty admissible binary word. There exists a unique triple $(v,a,b)$, where $v$ is an admissible binary word, empty or not, and $a$, $b$ are positive integers, such that $w=0\{1\}_{b-1}v\{0\}_{a-1}1$. For all integers $m\geqslant 1$ and $n\geqslant 1$, we have
\begin{equation}\label{E34}
\zeta(w)_{m-1,n-1}=\zeta(w)_{m,n}+n^{-a}\zeta(w^{\rm init})_{m,n}+m^{-b}\zeta(w^{\rm fin})_{m,n}+
n^{-a}m^{-b}\zeta(w^{\rm mid})_{m,n}
\end{equation}
where $w^{\rm init}=0\{1\}_{b-1}v$, $w^{\rm fin}=v\{0\}_{a-1}1$ and $w^{\rm mid}=v$.}
\end{prop}
By definition, $w$ starts by $0$ and ends by $1$. Any triple $(v,a,b)$, where $v$ is is an admissible binary word and $a$, $b$ are positive integers, such that $w=0\{1\}_{b-1}v\{0\}_{a-1}1$, is  necessarily 
obtained by the following procedure :\par
$a)$ one removes from $w$ the first bit $0$ and the last bit $1$ and one gets a word $w'$; \par
$b)$ the number of bits equal to $0$ by which $w'$ starts is $a-1$ and the number of bits equal to $1$ by which $w'$ ends is $b-1$;\par
$c)$ by removing these bits from $w'$, one gets $v$.\par
\noindent This shows the uniqueness of the triple $(a,b,v)$ as well as its existence, since the word $v$ obtained by the preceding procedure
either is empty, or starts by $0$ and ends by~$1$. By Proposition \ref{P3}, we have 
$$\zeta(w)_{m-1,n-1}=\zeta(w)_{m-1,n}+n^{-a}\zeta(w^{\rm init})_{m-1,n}$$
and by Proposition \ref{P4}, 
$${\zeta(w)_{m-1,n}=\zeta(w)_{m,n}+m^{-b}\zeta(w^{\rm fin})_{m,n},}\;\;$$
$${\zeta(w^{\rm init})_{m-1,n}=\zeta(w^{\rm init})_{m,n}+m^{-b}\zeta(w^{\rm mid})_{m,n}.}$$
Proposition \ref{P5} follows.
\bigskip

\section{An algorithm to compute  multiple zeta values}\label{S5}

Let $w=\varepsilon_1 \ldots \varepsilon_k$ be a non empty admissible binary word. Let ${\rm V}$ denote the set of non empty admissible subwords of $w,$ {\it i.e.~} of words
of the form $\varepsilon_i\ldots\varepsilon_j$, where $1\leqslant i\leqslant j\leqslant k$ and $\varepsilon_i=0$, $\varepsilon_j=1$. Let ${\rm N}$ be an integer.
We set $u_{\rm N}(v)=0$ for all $v\in{\rm V}$, and then compute inductively $u_n(v)$ for $v\in{\rm V}$, when $n$ is decreasing from ${\rm N}$ to $0$, by using the recurrence relation 
\begin{equation}\label{E35}
\hspace{0.5cm}
u_{n-1}(v)=u_n(v)+n^{-a(v)}u_n(v^{\rm init})+n^{-b(v)}u_n(v^{\rm fin})+n^{-a(v)-b(v)}u_n(v^{\rm mid}).
\end{equation}
In this formula, $v^{\rm init}$, $v^{\rm fin}$, $v^{\rm mid}$ are the initial, final  and  middle parts of $v$, as
defined in Theorem~ \ref{T5}, we have $a(v)=|v|-|v^{\rm init}|$, 
$b(v)=|v|-|v^{\rm fin}|$, and  in case 
$v^{\rm init}=0$, $v^{\rm fin}=1$ or $v^{\rm mid}=\varnothing$, the corresponding value  $u_n(0)$, $u_n(1)$ or $u_n(\varnothing)$ is taken to be $\zeta(0)_{n,n}$, $\zeta(1)_{n,n}$ or$\zeta(\varnothing)_{n,n}$, as defined by formula
(\ref{E14}), (\ref{E15}) or (\ref{E16}) respectively.

The real numbers $u_0(v)$ obtained when performing this algorithm are approximate values of the multiple zeta values $\zeta(v)$ for all $v\in{\rm V}$, and in particular for $w$. Errors are bounded as follows :\medskip

\begin{prop}\label{P6} The theoretical error $|\zeta(v)-u_0(v)|$ in the previous algorithm is bounded above by
 $2^{-2{\rm N}}({\rm N}+1)^2{\pi^2\over 6}$ for each $v\in{\rm V}$. If at each step of the algorithm
the right hand side of $(35)$ is computed  to an accuracy at most $ \alpha$, the total error (theoretical error plus rounding errors) is  bounded above by $2^{-2{\rm N}}({\rm N}+1)^2{\pi^2\over 6}+
{{\rm N}({\rm N}+1)(2{\rm N}+1)\over 6}\alpha$.
\end{prop}

We shall prove the second assertion. The first one is then obtained by taking $\alpha=0$. 
For $0\leqslant n\leqslant {\rm N}$, let $\beta_n$ denote the supremum, for $v\in{\rm V}$, of the numbers $|\zeta(v)_{n,n}-\underline{u}_n(v)|$, where $\underline{u}_n(v)$ is the computed value of 
$u_n(v)$. We have 
$$\beta_{\rm N}=\sup_{v\in{\rm V}} \zeta(v)_{{\rm N},{\rm N}}\leqslant 2^{-2\rm N}{\pi^2\over 6}$$ by Theorem \ref{T3},  and 
one deduces from formula $(35)$ that
$$\beta_{n-1}\leqslant (1+{1\over n}+{1\over n}+{1\over n^2})\beta_n+\alpha= \left({n+1\over n}\right)^2\beta_n+\alpha\qquad\quad\hbox{
for $1\leqslant n\leqslant {\rm N}$.}$$ It follows, by descending induction, that
$$(n+1)^2\beta_n\leqslant 2^{-2{\rm N}}({\rm N+1})^2{\pi^2\over 6}+\alpha \sum_{k=n+1}^{\rm N} k^2\qquad\hbox{
for $0\leqslant n\leqslant {\rm N}$.}$$
Proposition \ref{P6} is then obtained by taking $n=0$ in this inequality. \medskip

{\it Example}--- To compute $\zeta(w)$ with a precision of $10^{-1000}$, it suffices to perform
the algorithm with ${\rm N}=1673$, each step being computed to an accuracy of $10^{-1010}$.
\medskip

We see on this example that the accuracy of the computation has to be only slightly smaller than the required precision, and that the number ${\rm N}$ of required steps  is only slightly higher than ${\ln 10\over\ln 4}d$, where $10^{-d}$ is the required precision.

The algorithm presented here is especially efficient if one is interested in computing simultaneously all multiple zeta values up to a given weight $k$. It suffices  in that case to replace the preceding set ${\rm V}$ by the set of all non empty admissible binary words of weight at most $k$. 

I first implemented this algorithm in the language Python 2.7.5 on my personal computer. As an example, computing
the 127 multiple zeta values corresponding to admissible compositions of weight $\leqslant 8$ with 1000 exact decimal digits took 5 minutes and 9 seconds. With only 100 exact decimal digits, it took 5.8 seconds.
The same computations have been implemented by Henri Cohen in Pari/GP and in C , which have the advantage of being interpreted languages. They then take 0.9 seconds and 0.006 seconds respectively.
\bigskip
\section{Some examples}\label{S6}

\noindent {\it In weight 2}\par
As a first example, let us apply Theorem \ref{T5} to the word $w=01$, corresponding to the composition $(2)$. We have in  this case $w^{\rm init}=0$, $w^{\rm fin}=1$ and $w^{\rm mid}=\varnothing$, and for each integer $n\geqslant 1$, the recurrence relation expressed in terms of words
\begin{eqnarray}\label{E36}
\zeta(01)_{n-1,n-1}&=&\zeta(01)_{n,n}+n^{-1}\zeta(0)_{n,n}+n^{-1}\zeta(1)_{n,n}+
n^{-2}\zeta(\varnothing)_{n,n}
\nonumber\\&=&\zeta(01)_{n,n}+3n^{-2}\zeta(\varnothing)_{n,n}\,, 
\end{eqnarray}
or equivalently, in term of compositions,
\begin{equation}\label{E37}
\zeta(2)_{n-1,n-1}=\zeta(2)_{n,n}+3n^{-2}\zeta(\varnothing)_{n,n}\,,
\end{equation}
where $\zeta(\varnothing)_{n,n}=\left({2n\atop n}\right)^{-1}\!$.
We therefore have, for all integers $n\geqslant 0$,
\begin{equation}\label{E38}
\zeta(2)=\zeta(2)_{n,n}+3\sum_{m=1}^n m^{-2}\left({2m\atop m}\right)^{-1},\qquad 
\zeta(2)_{n,n}=3\sum_{m=n+1}^\infty m^{-2}\left({2m\atop m}\right)^{-1}.
\end{equation}

This yields in particular for $n=0$ the following formula, due to Euler :
\begin{equation}\label{E39}
\zeta(2)=3\sum_{m=1}^\infty m^{-2}\left({2m\atop m}\right)^{-1}.
\end{equation}

\noindent {\it In weight 4}\par
As a second example, let us apply Theorem \ref{T5} to the binary words $0001$, $0011$ and $0101$ corresponding respectively to the compositions
$(4)$, $(3,1)$ and $(2,2)$ :  The recurrence formulas obtained in the same way as in our first example can be expressed, after using the the duality relation $\zeta(2,1)_{n,n}=\zeta(3)_{n,n}$,  in the following matrix form :
\begin{equation}
\label{E40}
{\rm X}_{n-1}={\rm X}_n+{\rm A}{\rm Y}_n,
\end{equation}
where 
\begin{equation}
\label{E41}
{\rm X}_n=\pmatrix{\zeta(4)_{n,n}\cr\zeta(3,1)_{n,n}\cr\zeta(2,2)_{n,n}\cr}, \qquad {\rm A}=\pmatrix{1&0&2\cr 2&1&0\cr 0&2&1\cr}, \qquad
{\rm Y}_n=\pmatrix{n^{-1}\zeta(3)_{n,n}\cr n^{-2}\zeta(2)_{n,n}\cr n^{-4}\zeta(\varnothing)_{n,n}\cr}.
\end{equation}
In writing these relations, we have implicitly used the the duality relation\\ $\zeta(2,1)_{n,n}=\zeta(3)_{n,n}$. Since we have
\begin{equation}\label{E42}
\pmatrix{4&\!\!-2&1\cr}{\rm A}=\pmatrix{0&0&9\cr}, 
\end{equation}
the sequence $(u_n)_{n\in\bf N}$ defined by
\begin{equation} \label{E43}
u_n=\pmatrix{4&\!\!-2&1\cr}{\rm X}_n=4\zeta(4)_{n,n}-2\zeta(3,1)_{n,n}+\zeta(2,2)_{n,n}
\end{equation}
satisfies the recurrence relation 
\begin{equation}\label{E44}
u_{n-1}=u_n+9n^{-4}\left({2n\atop n}\right)^{-1}.
\end{equation}
Since it is known that $\zeta(3,1)={1\over 4}\zeta(4)$ and $\zeta(2,2)={3\over 4}\zeta(4)$, we have  $u_0={17\over 4}\zeta(4)$, and  therefore 
\begin{equation}
\label{E45}
\zeta(4)={36\over 17}\sum_{n=1}^\infty n^{-4}\left({2n\atop n}\right)^{-1}.
\end{equation}

\noindent {\it In higher  weight}\par
For any weight $k\geqslant 1$, we can as in our previous example consider the column vector ${\rm X}_n$, whose entries are the $(n,n)$-tails $\zeta({\bf a})_{n,n}$, where
$\bf a$ runs over a system of representatives up to duality of compositions of weight $k$. The corresponding recurrence relations
are again expressible in matrix form as
\begin{equation}
{\rm X}_{n-1}={\rm X}_n+{\rm A}{\rm Y}_n,
\end{equation}
where ${\rm Y}_n$ is the column vector whose entries are the numbers $n^{k_{\bf b-k}}\zeta({\bf b})_{n,n}$, where ${\bf b}$ runs  over a system of
representatives up to duality of compositions of weight $k_{\bf b}<k$, the last entry being  $n^{-k}\zeta(\varnothing)_{n,n}$. When $k=2k'$ is even, the  numbers of entries of ${\rm X}_n$ and ${\rm Y}_n$ are both equal to $2^{k'-2}(2^{k'-1}+1)$, so that ${\rm A}$ is a square matrix of this size.
As an example when $k=6$, we have
\bigskip
$${\rm X}_n=\!\!\pmatrix{\zeta(6)_{n,n}\cr\zeta(5,1)_{n,n}\cr\zeta(4,2)_{n,n}\cr\zeta(4,1,1)_{n,n}\cr
\zeta(3,3)_{n,n}\cr\zeta(3,2,1)_{n,n}\cr\zeta(3,1,2)_{n,n}\cr
\zeta(2,4)_{n,n}\cr\zeta(2,2,2)_{n,n}\cr\zeta(2,1,3)_{n,n}\cr}\!\!,\ \ {\rm A}=\!\!\pmatrix{1&0&0&0&0&0&0&0&0&2\cr
1&1&0&0&1&0&0&0&0&0\cr 0&0&1&0&1&0&0&1&0&0\cr 0&2&0&0&0&1&0&0&0&0\cr 0&0&0&1&0&0&0&1&1&0\cr
0&0&2&0&0&0&1&0&0&0\cr 0&0&0&1&0&1&0&1&0&0\cr 0&0&0&0&1&0&0&0&1&1\cr 0&0&0&0&0&0&2&0&1&0\cr
0&0&0&0&0&0&0&2&0&1\cr}\!\!,\ \ 
{\rm Y}_n =\!\!\pmatrix{n^{-1}\zeta(5)_{n,n}\cr n^{-1}\zeta(4,1)_{n,n}\cr n^{-1}\zeta(3,2)_{n,n}\cr n^{-1}\zeta(2,3)_{n,n}\cr n^{-2}\zeta(4)_{n,n}\cr
n^{-2}\zeta(3,1)_{n,n}\cr n^{-2}\zeta(2,2)_{n,n}\cr n^{-3}\zeta(3)_{n,n}\cr n^{-4}\zeta(2)_{n,n}\cr n^{-6}\zeta(\varnothing)_{n,n}\cr}\!\!.$$
\bigskip

However some new feature arises when $k=6$ : the last column of ${\rm A}$ is a linear combination of the other ones, and this forbids us to get a closed   formula
for $\sum_{n=1}^\infty n^{-6}\left({2n\atop n}\right)^{-1}$ as in the case of weight $4$. In fact the rank of ${\rm A}$ is $9$
and up to a scalar $${\rm L}=\pmatrix{2&\!\!-2&4&1&1&\!\!-2&\!\!-1&\!\!-2&1&\!\!-2\cr}$$
is the unique line vector such that ${\rm L}\rm A=0$. The sequence $(u_n)_{n\in\bf N}$ defined by ${\rm L}{\rm X}_n=(u_n)$ now satisfies the recurrence 
relation $u_{n-1}=u_n$ for $ n \geqslant 1$, hence is constant. Since it converges to~$0$, it is indentically $0$. In other words, the  linear combination of $(n,n)$-tails
\begin{eqnarray}
 2\zeta(6)_{n,n}-2\zeta(5,1)_{n,n}+4\zeta(4,2)_{n,n}+\zeta(4,1,1)_{n,n}+
\zeta(3,3)_{n,n}-2\zeta(3,2,1)_{n,n}\nonumber\\-\zeta(3,1,2)_{n,n}-2  
\zeta(2,4)_{n,n}+\zeta(2,2,2)_{n,n}-2\zeta(2,1,3)_{n,n} \qquad\qquad\qquad\qquad\;\;
\nonumber
\end{eqnarray}
vanishes for every~integer $n\in\bf N\,$!

Similar ${\bf Z}$-linear combinations of  compositions of weight k, for which the $(n,n)$-tails  vanish identically, occur for other weights $k$. We have computed by two different methods a lower bound for the rank $d_{k}$ of the $\bf{Z}$-module consisting of these linear combinations when $2 \leqslant k \leqslant 16$.  In~each case we found the same value, and we believe it to be the exact values of $d_{k}$.
\vspace{.4 cm}

\begin{table}[h]
\centering
\begin{tabular}{|c|c|c|c|c|c|c|c|c|c|c|c|c|c|c|c|}
\hline
$k$ & 2 & 3 & 4 & 5 & 6 & 7 & 8 & 9 & 10 & 11 & 12 & 13 & 14 & 15 & 16 \\
\hline
 $d_k\geqslant$& 0 & 0 & 0 & 0 & 1 & 0 & 4 & 2 & 14 & 15 & 52 & 78 & 200 & 350& 789 \\
\hline
\end{tabular} 
\end{table}

We do not have at present a theoretical understanding of these numerical observations and we hope to be able to address them in a forthcoming paper. We were not able to find any weight $k\geqslant 5$ for which $\sum_{n=1}^{\infty} n^{-k} \left( 2n \atop n \right)^{-1}$ could be expressed as a ${\bf Z}$-linear combination of multiple zeta values.

\bigskip

\section{Iterating the recurrence relations for double tails}\label{S7}

Our purpose in this section is to write down the formula obtained after ${\rm N}$ iterations of formula (\ref{E34}).
Before stating our result, we shall need to introduce some notations.\medskip

\subsection{Two notations}\label{S7.1}

For any non empty composition ${\bf a}=(a_1,\ldots,a_r)$ and any integers $p$, $q$ such that\\ $0\leqslant p< q$,
 we define a real number $\varphi_{p,q}(\bf a)$ by the finite sum
 \begin{equation}\label{E47}
\varphi_{p,q}({\bf a})=q^{-a_1}\sum_{q>n_2>\ldots>n_r>p} n_2^{-a_2}\ldots n_r^{-a_r},
\end{equation}
where the right hand side is considered to be equal to $q^{-a_1}$ when $r=1$.

For any  composition ${\bf a}=(a_1,\ldots,a_r)$ and any integers $p$, $q$ such that $0\leqslant p\leqslant q$, we define a real number $\zeta_{]p,q]}(\bf a)$ by  the finite sum
\begin{equation}\label{E48}
\zeta_{]p,q]}({\bf a})=\sum_{q\geqslant n_1>n_2>\ldots>n_r>p} n_1^{-a_1} n_2^{-a_2}\ldots n_r^{-a_r},
\end{equation}
where the right hand side is considered to be equal to $1$ when $r=0$.

When $\bf a$ is non empty, we clearly have
\begin{equation}\label{E49}
\zeta_{]p,q]}({\bf a})=\sum_{p<q'\leqslant q}\varphi_{p,q'}({\bf a}),
\end{equation}
and if $q>p$,
\begin{equation} \label{E50}
\varphi_{p,q}({\bf a})=q^{-a_1}\zeta_{]p,q-1]}(a_2,\ldots,a_r).
\end{equation}
\par

We shall need the following estimations :\medskip
\begin{prop}---
\label{P7}  Assume $0\leqslant p< q$. We have $\varphi_{p,q}({\bf a})\leqslant 1$ when $\bf a$ is non empty, and $\zeta_{]p,q]}({\bf a})\leqslant q-p$ for all $\bf a$.
\end{prop} 

The second assertion is true when $\bf a$ is empty. If the second assertion is true
for compositions of a given depth $r$, the first one is true for compositions of depth $r+1$ by 
formula (\ref{E50}), and then the second one is  true for compositions of depth $r+1$ by 
formula~(\ref{E49}). This proves the proposition by induction. \bigskip

\subsection{Statement of the result}\label{S7.2}

In this section, $w=\varepsilon_1\ldots\varepsilon_k$ denotes a non empty admissible binary word. 

 For any pair $s=(i,j)$ of integers such that
$0\leqslant i\leqslant j\leqslant k$, we denote by $w_{s}$ the binary word $\varepsilon_{i+1}\ldots\varepsilon_j$, by ${\bf a}_s$ the composition corresponding to the binary word $\varepsilon_{j+1}\ldots\varepsilon_n$ and by ${\bf b}_s$ the composition corresponding to the dual of the binary word $\varepsilon_{1}\ldots\varepsilon_i$. Hence  $w$ is equal to $\overline{{\bf w}({\bf b}_s)}w_s{\bf w}({\bf a}_s)$.

Let ${\rm S}$ (resp. ${\rm S}_0$; resp. ${\rm S}_1$) denote the set of pairs  of integers $s=(i,j)$ such that
$0\leqslant i\leqslant j\leqslant k$ for which $w_s\in{}_0{\rm W}_{1}$ (resp. $w_s=0$; resp.  $w_s=1$). 
Finally, let ${\rm S}_\varnothing$ denote the set of pairs $(i,i)$, with $0< i< k$.
With these notations :\medskip

\begin{thm}--- \label{T7}For all $m\geqslant 0$, $n\geqslant 0$ and ${\rm N}\geqslant 1$, we have
\begin{eqnarray}
\label{E51}
\zeta(w)_{m,n}&=&\sum_{s\in{\rm S}}\zeta_{]m,m+{\rm N}]}({\bf b}_s)\zeta_{]n,n+{\rm N}]}({\bf a}_s)
\zeta(w_s)_{m+{\rm N},n+{\rm N}}\;\;\;\;\;\;\;\;\;\;\;\;\;\;
\nonumber \\&+&\sum_{\rm M=1}^{\rm N}\sum_{s\in{\rm S}_0}\zeta_{]m,m+{\rm M}-1]}({\bf b}_s)\varphi_{n,n+{\rm M}}({{\bf a}_s})
\zeta(0)_{m+{\rm M},n+{\rm M}}
\nonumber\\
&+&\sum_{{\rm M}=1}^{\rm N}\sum_{s\in{\rm S}_1}\varphi_{m,m+{\rm M}}({{\bf b}_s})\zeta_{]n,n+{\rm M}-1]}({{\bf a}_s})
\zeta(1)_{m+{\rm M},n+{\rm M}}
\\&+&\sum_{{\rm M}=1}^{\rm N}\sum_{s\in{\rm S}_\varnothing}\varphi_{m,m+{\rm M}}({\bf b}_s)\varphi_{n,n+{\rm M}}({\bf a}_s)
\zeta(\varnothing)_{m+{\rm M},n+{\rm M}}.\nonumber
\end{eqnarray}
\end{thm}

\subsection{Proof of Theorem 7}\label{7.3}

For each ${\rm N}\geqslant 0$, let 
$\rm A_{\rm N}=\sum_{s\in{\rm S}}\zeta_{]m,m+{\rm N}]}({\bf b}_s)\zeta_{]n,n+{\rm N}]}({\bf a}_s)
\zeta(w_s)_{m+{\rm N},n+{\rm N}}.$

When ${\rm N}=0$, the only non zero term in the previous sum is the one corresponding to $s=(0,k)$, for which 
${\bf a}_s={\bf b}_s=\varnothing$ and $w_s=w$. We hence have ${\rm A}_0=\zeta(w)_{m,n}$.

In order to prove Theorem \ref{T7} by induction, it therefore suffices to prove that

\begin{eqnarray} \label{E52}
{\rm A}_{{\rm N}-1}={\rm A}_{\rm N} 
&+&\sum_{s\in{\rm S}_0}\zeta_{]m,m+{\rm N}-1]}({\bf b}_s)\varphi_{n,n+{\rm N}}({\bf a}_s)
\zeta(0)_{m+{\rm N},n+{\rm N}}
\nonumber \\ 
&+&\sum_{s\in{\rm S}_1}\varphi_{m,m+{\rm N}}({\bf b}_s)\zeta_{]n,n+{\rm N}-1]}({\bf a}_s)
\zeta(1)_{m+{\rm N},n+{\rm N}}
 \\
& +&\sum_{s\in{\rm S}_\varnothing}\varphi_{m,m+{\rm N}}({\bf b}_s)\varphi_{n,n+{\rm N}}({\bf a}_s)
\zeta(\varnothing)_{m+{\rm N},n+{\rm N}}\,\;\nonumber
\end{eqnarray}
holds for ${\rm N}\geqslant 1$. 

From now on, we assume ${\rm N}\geqslant 1$ and denote for simplicity $m+{\rm N}$ and $n+{\rm N}$
by $\underline{m}$ and $\underline{n}$. 
For each $s=(i,j)\in{\rm S}$, there exists 
a unique pair of integers $(i',j')$ such that $i\leqslant i'\leqslant j'\leqslant j$ 
and such that the initial, final and middle parts of $w_s$ are $w_{s^{\rm init}}$, $w_{s^{\rm fin}}$
and $w_{s^{\rm mid}}$, where
 $s^{\rm init}=(i,j')$,  $s^{\rm fin}=(i',j)$ and  $s^{\rm mid}=(i',j')$. Proposition 5 then yields
$$\zeta(w_s)_{\underline{m}-1,\underline{n}-1}=\zeta(w_s)_{\underline{m},\underline{n}}
+\underline{n}^{j'-j}\zeta(w_{s^{\rm init}})_{\underline{m},\underline{n}}+\underline{m}^{i-i'}\zeta(w_{s^{\rm fin}})_{\underline{m},\underline{n}}
+\underline{m}^{i-i'}\underline{n}^{j'-j}\zeta(w_{s^{\rm mid}})_{\underline{m},\underline{n}}.$$
We  have ${\bf a}_{s^{\rm fin}}={\bf a}_{s}$, ${\bf b}_{s^{\rm init}}={\bf b}_{s}$,
${{\bf a}_{s^{\rm mid}}}={\bf a}_{s^{\rm init}}$, ${\bf b}_{s^{\rm mid}}={\bf }b_{s^{\rm fin}}$ and by formula (50)
$$\underline{n}^{j'-j}\zeta_{]n,\underline{n}-1]}({\bf a}_s)=\varphi_{n,\underline{n}}({\bf a}_{s^{\rm init}}),\qquad
\underline{m}^{i-i'}\zeta_{]m,\underline{m}-1]}({\bf b}_s)=\varphi_{m,\underline{m}}({{\bf b}_{s^{\rm fin}}}).$$
Combining these equalities with the previous one, we see that 
$$\zeta_{]m,\underline{m}-1]}({\bf b}_s)\zeta_{]n,\underline{n}-1]}({\bf a}_s)
\zeta(w_s)_{\underline{m}-1,\underline{n}-1}$$ is equal to
\begin{eqnarray}
& &\zeta_{]m,\underline{m}-1]}({\bf b}_s)\zeta_{]n,\underline{n}-1]}({\bf a}_s)
\zeta(w_s)_{\underline{m},\underline{n}}
+\zeta_{]m,\underline{m}-1]}({\bf b}_{s^{\rm init}})\varphi_{n,\underline{n}}({\bf a}_{s^{\rm init}})
\zeta(w_{s^{\rm  init}})_{\underline{m},\underline{n}}
\quad\;\nonumber\\&+&\varphi_{m,\underline{m}}({\bf b}_{s^{\rm fin}})\zeta_{]n,\underline{n}-1]}({\bf a}_{s^{\rm fin}})
\zeta(w_{s^{\rm fin}})_{\underline{m},\underline{n}}
+\varphi_{m,\underline{m}}({\bf b}_{s^{\rm mid}})\varphi_{n,\underline{n}}({\bf a}_{s^{\rm mid}})
\zeta(w_{s^{\rm mid}})_{\underline{m},\underline{n}}.\nonumber
\end{eqnarray}

Note that for any $s\in{\rm S}$, the pairs $s^{\rm init}$, $s^{\rm fin}$ and $s^{\rm mid}$ belong to
$\rm S\cup{\rm S}_0\cup{\rm S}_1\cup{\rm S}_{\varnothing}$, and ${\rm A}_{{\rm N}-1}$ is the sum, for $s$ in ${\rm S}$, of the previous expressions. We shall get a new formula for ${\rm A}_{{\rm N}-1}$
by adding all terms involving $\zeta(w_t)_{\underline{m},\underline{n}}$ for a given 
$t$ in ${\rm S}\cup{\rm S}_0\cup{\rm S}_1\cup{\rm S}_{\varnothing}$, and then summing over $t$.
We distinguish several cases :\smallskip
\\
\noindent $a)$ {\it Case where $t=(i,j)\in{\rm S}$}\par
\noindent --- If $i=0$ and $j=k$, there is  a single term containing $\zeta(w_t)_{\underline{m},\underline{n}}$; it corresponds to $s=t$. Its coefficient is 
$$\zeta_{]m,\underline{m}-1]}({\bf b}_t)\zeta_{]n,\underline{n}-1]}({\bf a}_t)=1=
\zeta_{]m,\underline{m}]}({\bf b}_t)\zeta_{]n,\underline{n}]}({\bf a}_t).$$
---  If $i=0$ and $j<k$, there are two terms containing $\zeta(w_t)_{\underline{m},\underline{n}}$, one for which $t=s$ and one for which $t=s^{\rm init}$. The sum of their coefficients is 
$$\zeta_{]n,\underline{n}-1]}({\bf a}_t)+\varphi_{n,\underline{n}}({\bf a}_t)=\zeta_{]n,\underline{n}]}({\bf a}_t)=
\zeta_{]m,\underline{m}]}({\bf b}_t)\zeta_{]n,\underline{n}]}({\bf a}_t).$$
---  If $i>0$ and $j=k$, there are two terms containing $\zeta(w_t)_{\underline{m},\underline{n}}$, one for which $t=s$ and one for which $t=s^{\rm fin}$. The sum of their coefficients is 
$$\zeta_{]m,\underline{m}-1]}({\bf b}_t)+\varphi_{m,\underline{m}}({\bf b}_t)
=\zeta_{]m,\underline{m}]}({\bf b}_t)=
\zeta_{]m,\underline{m}]}({\bf b}_t)\zeta_{]n,\underline{n}]}({\bf a}_t).$$
---  If $i>0$ and $j<k$, there are four terms containing $\zeta(w_t)_{\underline{m},\underline{n}}$,  for which $t$ is equal to $s$, $s^{\rm init}$,  $s^{\rm fin}$ and $s^{\rm mid}$
respectively. The sum of their coefficients is 
$${(\zeta_{]m,\underline{m}-1]}({\bf b}_t)+\varphi_{m,\underline{m}}({\bf b}_t))
(\zeta_{]n,\underline{n}-1]}({\bf a}_t)+\varphi_{n,\underline{n}}({\bf a}_t))
=\zeta_{]m,\underline{m}]}({\bf b}_t)\zeta_{]n,\underline{n}]}({\bf a}_t).}$$
\noindent $b)$ {\it Case where $t\in{\rm S}_0$}\par
In this case there is only one term containing $\zeta(w_t)_{\underline{m},\underline{n}}$; it corresponds to an $s$ for which $t=s^{\rm init}$. Its coefficient is $\zeta_{]m,\underline{m}-1]}({\bf b}_t)\varphi_{n,\underline{n}}({\bf a}_t)$.\smallskip
\\
\noindent $c)$ {\it Case where $t\in{\rm S}_1$}\par
In this case there is only one term containing $\zeta(w_t)_{\underline{m},\underline{n}}$; it corresponds to an $s$ for which $t=s^{\rm fin}$. Its coefficient is $\varphi_{m,\underline{m}}({\bf b}_t)
\zeta_{]n,\underline{n}-1]}({\bf a}_t)$.\smallskip
\\
\noindent $d)$ {\it Case where $t\in{\rm S}_\varnothing$}\par
In this case there is only one term containing $\zeta(w_t)_{\underline{m},\underline{n}}$; it corresponds to an $s$ for which $t=s^{\rm mid}$. Its coefficient is $\varphi_{m,\underline{m}}({\bf b}_t)
\varphi_{n,\underline{n}}({\bf a}_t)$.\smallskip

Summing up all the terms, we get
\begin{eqnarray} 
{\rm A}_{{\rm N}-1}&=&\sum_{t\in{\rm S}}\zeta_{]m,\underline{m}]}({\bf b}_t)\zeta_{]n,\underline{n}]}({\bf a}_t)
\zeta(w_t)_{\underline{m},\underline{n}}+\sum_{t\in{\rm S}_0}\zeta_{]m,\underline{m}-1]}({\bf b}_t)\varphi_{n,\underline{n}}({\bf a}_t)
\zeta(0)_{\underline{m},\underline{n}}
\nonumber \\&+&\sum_{t\in{\rm S}_1}\varphi_{m,\underline{m}}({\bf b}_t)\zeta_{n,\underline{n}-1}({\bf a}_t)
\zeta(1)_{\underline{m},\underline{n}}
+\sum_{t\in{\rm S}_\varnothing}\varphi_{m,\underline{m}}({\bf b}_t)\varphi_{n,\underline{n}}({\bf a}_t)
\zeta(\varnothing)_{\underline{m},\underline{n}},\;\;\nonumber
\end{eqnarray}
which is exactly the relation (\ref{E52}) we wanted to establish. \bigskip

\subsection{Reformulation of Theorem 7}\label{S7.4}

We keep the same notations as in section \ref{S7.2}.  In particular $w=\varepsilon_1\ldots\varepsilon_k$ is a non empty admissible binary word.
 
For $1\leqslant i\leqslant k-1$, let ${\bf a}_i$ denote the composition corresponding to the binary
word $\varepsilon_{i+1}\ldots\varepsilon _k$ and ${\bf b}_i$ the composition corresponding to the dual of the binary word $\varepsilon_{1}\ldots\varepsilon _i$. In other word, we have 
${\bf a}_i={\bf a}_s$ and ${\bf b}_i={\bf b}_s$, where $s$ is the element $(i,i)$ of ${\rm S}_\varnothing$.

 For any pair $(\varepsilon,\varepsilon')$ of bits and any positive integers $m$ and $n$, we define a number
 $\lambda_{m,n}(\varepsilon,\varepsilon ')$ by
 \begin{equation}
 \label{E53}
 \lambda_{m,n}(\varepsilon,\varepsilon ')=1+\cases{{m\over n}& if $\varepsilon=0$\cr
 0& if $\varepsilon=1$\cr}+\cases{{n\over m}& if $\varepsilon'=1$\cr
 0& if $\varepsilon'=0$.\cr}
  \end{equation}

 Theorem \ref{T7} can the be restated as follows :\bigskip
 
 \begin{thm}--- \label{T8}For all $m\geqslant 0$, $n\geqslant 0$ and ${\rm N}\geqslant 1$, we have
 \begin{equation}
  \label{E54}
 \hspace{0.3cm} \begin{array}{rrl} \zeta(w)_{m,n} & = & \Large{\sum_{s\in {\rm  S} } \zeta_{]m,m+{\rm N}]}({{\bf b}_s})\zeta_{]n,n+{\rm N}]}({\bf a}_s)
  \zeta(w_s)_{m+{\rm N},n+{\rm N}}} 
  \\[1em]  & + & \Large{\sum_{{\rm M}=1}^{{\rm N}} \left({m+n+2 {\rm M} \atop m+{\rm  M} }  \right)^{-1} \sum_{i=1}^{k-1} \varphi_{m,m+{\rm M}}({\bf b}_i) \varphi_{n,n+{\rm M}}({{\bf a}_i}) \lambda_{m+{\rm M},n+{\rm M}}(\varepsilon_i,\varepsilon_{i+1})}\end{array}
\end{equation}
\end{thm}
Theorem \ref{T8} follows from Theorem \ref{T7} and the following observations :\par
$a)$ The elements of ${\rm S}_0$ are the pairs $(i-1,i)$, where $1\leqslant i\leqslant k-1$ and 
$\varepsilon_i=0$. For such a pair $s$, we have ${\bf a}_i={\bf a}_s$, ${\bf b}_i=(1,{\bf b}_s)$, hence by formula (\ref{E50})
$$\zeta_{]m,m+{\rm M}-1]}({\bf b}_s)=(m+{\rm M})\varphi_{m,m+{\rm M}}({\bf b}_i),\qquad
\varphi_{n,n+{\rm M}}({\bf a}_s)=\varphi_{n,n+{\rm M}}({\bf a}_i).$$\par
$b)$ The elements of ${\rm S}_1$ are the pairs $(i,i+1)$, where $1\leqslant i\leqslant k-1$ and 
$\varepsilon_{i+1}=1$. For such a pair $s$,  we have ${\bf a}_i=(1,{\bf a}_s)$, ${\bf b}_i={\bf b}_s$, hence by formula (\ref{E50})
$$\varphi_{m,m+{\rm M}}({\bf b}_s)=\varphi_{m,m+{\rm M}}({\bf b}_i),\qquad
\zeta_{]n,n+{\rm M}-1]}({\bf a}_s)=(n+{\rm M})\varphi_{n,n+{\rm M}}({\bf a}_i).$$\par
$c)$ The elements of ${\rm S}_\varnothing$ are the pairs $(i,i)$, where $1\leqslant i\leqslant k-1$. For such a pair $s$, we have 
$$\varphi_{m,m+{\rm M}}({\bf b}_s)=\varphi_{m,m+{\rm M}}({\bf b}_i),\qquad
\varphi_{n,n+{\rm M}}({\bf a}_s)=\varphi_{n,n+{\rm M}}({\bf a}_i).$$\par
$d)$ We have $$\zeta(\varnothing)_{m+{\rm M},n+{\rm M}}=\left({m+n+2{\rm M}\atop m+{\rm M}}\right)^{-1}$$ and 
$$\zeta(0)_{m+{\rm M},n+{\rm M}}={\zeta(\varnothing)_{m+{\rm M},n+{\rm M}}\over n+{\rm M}}\cvirg\qquad
\zeta(1)_{m+{\rm M},n+{\rm M}}={\zeta(\varnothing)_{m+{\rm M},n+{\rm M}}\over m+{\rm M}}\cdot$$ \bigskip

\subsection{Proof of Theorem 6}\label{S7.5}

We keep the same notations as in section \ref{S7.4}.  In particular $w=\varepsilon_1\ldots\varepsilon_k$ is a non empty admissible binary word.
\medskip

 \begin{thm}--- \label{T9}For all $m\geqslant 0$, $n\geqslant 0$ and ${\rm N}\geqslant 1$, we have
 \begin{equation}
 \label{E55}\hspace{1cm}
\zeta(w)_{m,n}=\sum_{{\rm M}=1}^{+\infty}\left({m+n+2{\rm M}\atop m+{\rm M}}\right)^{-1}\sum_{i=1}^{k-1}\varphi_{m,m+{\rm} M}({\bf b}_i)\varphi_{n,n+{\rm M}}({\bf a}_i)\lambda_{m+{\rm M},n+{\rm M}}(\varepsilon_i,\varepsilon_{i+1}).
 \end{equation}
\end{thm}

We deduce this theorem from Theorem \ref{T8} by letting ${\rm N}$ tend to $+\infty$. It suffices to check that, for 
any given $s\in{\rm S}$, the product
$\zeta_{]m,m+\rm N]}({\bf b}_s)\zeta_{]n,n+{\rm N}]}({\bf a}_s)
\zeta(w_s)_{m+{\rm N},n+{\rm N}}$ tends to $0$ when ${\rm N}$ tends to $+\infty$. But this folllows form the fact that
$$0\leqslant \zeta_{]m,m+{\rm N}]}({\bf b}_s)\leqslant {\rm N},\qquad
0\leqslant\zeta_{]n,n+{\rm N}]}({\bf a}_s)\leqslant {\rm N},$$
for any ${\rm N}\geqslant 1$ by proposition \ref{P7}, and
$$0\leqslant \zeta(w_s)_{m+{\rm N},n+{\rm N}}\leqslant {(m+{\rm N})^{m+{\rm N}} (n+{\rm N})^{n+{\rm N}}\over
 (m+n+2{\rm N})^{m+n+2{\rm N}}}\zeta(w_s)$$
 for any ${\rm N}\geqslant 0$ by Theorem \ref{T3}.\medskip

Theorem \ref{T6} is the particular case of Theorem \ref{T9} where $m=n=0$.

\bigskip
\section{Another algorithm to compute  multiple zeta values}\label{S8}

Let ${\bf a}$ be a non empty admissible composition. Let $\varepsilon_1\ldots\varepsilon_k$ denote the corresponding binary word. Theorem 6 states that
\begin{equation}\label{E56}
\zeta({\bf a})=\sum_{m=1}^\infty\psi_m({\bf a})\left({2m\atop m}\right)^{-1}
\end{equation}
with
\begin{equation}\label{E57}
\psi_m({\bf a})=\sum_{i=1}^{k-1}\lambda(\varepsilon_i,\varepsilon_{i+1})\varphi_m({\bf a}_i)\varphi_m({\bf b}_i),
\end{equation}
where ${\bf a}_i$ and ${\bf b}_i$ denote the compositions corresponding to the binary words
$\varepsilon_{i+1}\ldots\varepsilon_k$ and $\overline{\varepsilon_i}\ldots\overline{\varepsilon_1}$, 
and where the functions $\varphi_m$ and $\lambda$ are defined by formulas (\ref{E23}) and (\ref{E24}).

Computing the partial sum 
\begin{equation}\label{E58}
\sum_{m=1}^{\rm N}\psi_m({\bf  a}) \left({2m\atop m}\right)^{-1}
\end{equation}
for some (large) integer ${\rm N}\geqslant 1$ yields an approximate value of $\zeta({\bf a})$. Since the real numbers $\varphi_m({\bf a}_i)$ and $\varphi_m({\bf b}_i)$ belong to $[0,1]$ by Proposition \ref{P7}, and $\lambda(\varepsilon_i,\varepsilon_{i+1})$ is an integer equal to $1$, $2$ or $3$, we have $0\leqslant \psi_m({\bf a})\leqslant 3(k-1)$. Therefore the theoretical error done by replacing the infinite sum by the partial sum up to ${\rm N}$ is bounded by 
$$3(k-1)\sum_{m={\rm N}+1}^\infty \left({2m\atop m}\right)^{-1}.$$

We have 
$$\left({2m\atop m}\right)^{-1}\leqslant m2^{1-2m}\leqslant (3m-4)2^{1-2m}
=2(m-1)2^{-2(m-1)}-2m2^{-2m}$$
for $m\geqslant 2$, hence the sum $\sum_{m={\rm N}+1}^\infty \left({2m\atop m}\right)^{-1}$ is bounded by
$2^{1-2{\rm N}}{\rm N}$ for ${\rm N}\geqslant 1$ and the theoretical error in our algorithm is bounded by $6(k-1)2^{-2{\rm N}}{\rm N}$.

The real numbers $\varphi_m({\bf a}_i)$, where $1\leqslant m\leqslant {\rm N}$ and $1\leqslant i\leqslant k-1$, are easily computed by induction on $m$, and for a given $m$, by descending induction on $i$. This is
done for example by using the following relations :
$$\varphi_m({\bf a}_i)=\cases{m^{-1}&if $i=k-1$,\cr
m^{-1}\varphi_m({\bf a}_{i+1})&if $\varepsilon_{i+1}=0$,\cr
0&if $m=1$, $\varepsilon_{i+1}=1$ and $i\not=k-1$,\cr
(1-m^{-1})\varphi_{m-1}({\bf a}_i)+m^{-1}\varphi_{m-1}({\bf a}_{i+1})&if $m\geqslant 2$, $\varepsilon_{i+1}=1$ and $i\not=k-1$.}$$

This algorithm, implemented in Python 2.7.5 on my personal computer took for example  10.1 seconds to compute $\zeta(2,1,3,2)$ with 1000 exact decimal places and 0.24 seconds with  only 100 exact decimal places. Implemented later by Henri Cohen in C it took 0.08 and 0.001 seconds respectively.

Another algorithm previously used to compute $\zeta({\bf a})$,  and for example implemented on the site EZ-face quoted in the introduction, is based on the following identity (a version of Chasles formula
for iterated integrals):
\begin{equation}\label{E59}
\zeta({\bf a})=\sum_{i=0}^k{\rm Li}_{{\bf a}_i}({1\over2}){\rm Li}_{{\bf b}_i}({1\over2})
\end{equation}
where ${\bf a}_i$ and ${\bf b}_i$ are defined as before and
$${\rm Li}_{(a_1,\ldots,a_r)}(z)=
\sum_{n_1>\ldots>n_r>0}{z^{n_1}\over n_1^{a_1}\ldots n_r^{a_r}}\cdot$$
for any composition $(a_1,\ldots,a_r)$ and any $z\in{\bf C}$ such that $|z|<1$.

Formula (\ref{E59}) on which the latter algorithm is based shows some structural similarities with formula  (\ref{E57}); we do not really understand their origin. However the first algorithm requires  only half as many steps as the second one to achieve the same precision, and implemented on the same computer, it takes roughly a third of the time. Moreover, the second algorithm looks somewhat artificial, because it involves splitting the iterated integral from $0$ to $1$ at ${1\over 2}$, whereas one could also choose any other intermediate element between $0$ and $1$.


\bigskip
{\bf Acknowledgement :}  I am deeply grateful to Professor Oesterl\'{e} for guiding me through this work.

\Addresses
\end{document}